\input amstex
\documentstyle{amsppt}
\document

\magnification 1100

\def\gen{{\frak{g}}}

\def\ten{\frak{t}}

\def\a{{\alpha}}

\def\eps{{\varepsilon}}

\def\1b{{\bold 1}}

\def\hb{{\bold h}}

\def\vb{{\bold v}}
\def\wb{{\bold w}}
\def\xb{{\bold x}}

\def\Ub{{\bold U}}

\def\Yb{{\bold Y}}

\def\vbt{{\bold {\tilde v}}}

\def\L{{\roman L}}

\def\K{{\roman K}}

\def\Hom{\text{Hom}\,}

\def\Res{\text{res}}

\def\Ker{\text{Ker}\,}

\def\rk{\text{rk\,}}

\def\Gr{\text{Gr}}

\def\Lie{\text{Lie}}
\def\GL{{\text{GL}}}

\def\CC{{\Bbb C}}

\def\NN{{\Bbb N}}
\def\PP{{\Bbb P}}

\def\ZZ{{\Bbb Z}}

\def\Ec{{\Cal E}}
\def\Fc{{\Cal F}}

\def\Lc{{\Cal L}}

\def\Oc{{\Cal O}}

\def\Vc{{\Cal V}}
\def\Wc{{\Cal W}}

\def\and{{\quad\text{and}\quad}}

\def\ds{\displaystyle}

\def\qed{\hfill $\sqcap \hskip-6.5pt \sqcup$}        
\overfullrule=0pt                                    

\def\u1{{\underline 1}}

\def\oh{{\overline h}}

\newdimen\Squaresize\Squaresize=14pt
\newdimen\Thickness\Thickness=0.5pt
\def\Square#1{\hbox{\vrule width\Thickness
	      \vbox to \Squaresize{\hrule height \Thickness\vss
	      \hbox to \Squaresize{\hss#1\hss}
	      \vss\hrule height\Thickness}
	      \unskip\vrule width \Thickness}
	      \kern-\Thickness}
\def\Vsquare#1{\vbox{\Square{$#1$}}\kern-\Thickness}

\noindent{\bf QUIVER VARIETIES AND YANGIANS}

\vskip2mm

\noindent MICHELA VARAGNOLO \footnote{
The author is partially supported by EEC grant
no. ERB FMRX-CT97-0100.}

\vskip3mm

\noindent {\it D\'epartement de Math\'ematique, Universit\'e de Cergy-Pontoise, 2 av. A. Chauvin,
95302 Cergy-Pontoise cedex. e-mail: michela.varagnolo\@math.u-cergy.fr}

\vskip3mm



\vskip1mm

\noindent{\bf   Abstract.}
We prove a conjecture of Nakajima ([5], for type $A$ it was announced in [2])
giving a geometric realization,
via quiver varieties, of the Yangian of type $ADE$ (and more in general of the Yangian
associated to every symmetric Kac-Moody Lie algebra).
As a corollary we get that the finite dimensional representation theory of the quantized affine
algebra and that of the Yangian coincide.

\vskip5mm

\noindent {\bf 1. The algebra $\Yb_\hbar(\L{\frak g})$.}

\noindent Let ${\frak g}$ be a simple, simply laced, complex Lie algebra over $\CC$
with Cartan matrix $A=(a_{kl})_{k,l\in I}.$
Denote by $L\gen=\gen[t,t^{- 1}]$ the loop Lie algebra of $\gen$. 
The Yangian  $\Yb_\hbar(\L{\frak g})$ is  the associative algebra, free
over $\CC[\hbar],$ generated
by $\xb^\pm_{k,r},\,\hb_{k,r}$  $(k\in I,\, r\in\NN)$
with the following defining relations
$$[\hb_{k,r},\hb_{l,s}]=0,\quad
[\hb_{k,0},\xb^\pm_{l,s}]=\pm a_{kl}\xb^\pm _{l,s},\leqno(1.1) $$
$$2[\hb_{k,r+1},\xb^\pm_{l,s}]-2[\hb_{k,r},\xb^\pm_{l,s+1}]=\pm\hbar
a_{kl}(\hb_{k,r}\xb^\pm_{l,s}+\xb^\pm_{l,s}\hb_{k,r}),\leqno(1.2)$$
$$[\xb^+_{k,r},\xb^-_{l,s}]=\delta_{kl}\hb_{k,r+s},\leqno(1.3)$$
$$2[\xb^\pm_{k,r+1},\xb^\pm_{l,s}]-2[\xb^\pm_{k,r},\xb^\pm_{l,s+1}]=
\pm\hbar a_{kl}(\xb^\pm_{k,r}\xb^\pm_{l,s}+\xb^\pm_{l,s}\xb^\pm_{k,r}),
\leqno(1.4)$$
$$\sum_{w\in S_m}[\xb^\pm_{k,r_{w(1)}},[\xb^\pm_{k,r_{w(2)}},...,[\xb^\pm_
{k,r_{w(m)}},\xb^\pm_{l,s}]...]]=0,\quad k\not= l\leqno(1.5)$$
\noindent for all sequences of non-negative integers $r_1,...,r_m$, where
$m=1-a_{kl}.$ 

\noindent Set 
$$[n]={{q^n-q^{-n}}\over{q-q^{-1}}} \qquad\forall n\in\ZZ.$$

\vskip2mm

\noindent{\bf 2. Quiver varieties.} 

\noindent Let $I$ (resp. $E$) be the set of vertices (resp. edges) of a finite
graph $(I,E)$ with no edge loops. For $k,l\in I$ let $n_{kl}$ be
the number of edges joining $k$ and $l$. Put $a_{kl}=2\delta_{kl}-n_{kl}$.
The map $(I,E)\mapsto A=(a_{kl})_{k,l\in I}$ is a bijection from the set of
finite graphs with no loops onto the set of symmetric  generalized
Cartan matrices. Let $\a_k$ and $\omega_k$, $k\in I$,
be the simple roots and fundamental weights 
of the symmetric Kac-Moody algebra corresponding to $A.$
Let $H$ be the set of edges of $(I,E)$ together 
with an orientation. For $h\in H$ let $h'\in I$ (resp. $h''\in I$) 
the incoming (resp. the outcoming) vertex of $h$. If $h\in H$ we denote
by $\oh\in H$ the same edge with opposite orientation.
Take two collection of finite dimensional complex  vector spaces 
$V=(V_k)_{k\in I},\, W=(W_k)_{k\in I}.$ 
Let us fix once for all the following convention :
the dimension of the graded vector space $V$
is identified with the element $\vb=\sum_{k\in I}v_k\a_k$ in the root lattice
(where $v_k$ is the dimension of $V_k$). 
Similarly the dimension of $W$ is 
identified with the weight $\wb=\sum_kw_k\omega_k$ 
(where $w_k$ is the dimension of $W_k$). Set
$$M(\vb,\wb)= \bigoplus_{h\in H}\Hom(V_{h''},V_{h'})
\oplus 
\bigoplus_{k\in I}\Hom(W_k,V_k)
\oplus 
\bigoplus_{k\in I}\Hom(V_k,W_k).$$
The group $G_\vb=\prod_k\GL(V_k)$ acts on $M(\vb,\wb)$ by 
$g\cdot(B,i,j)=(gBg^{-1},gi,jg^{-1}).$
We denote by $B_h$ the component of the element $B$
in $\Hom(V_{h''},V_{h'})$. Let us consider the map
$$\mu_{\vb,\wb}\,:\, M(\vb,\wb)\to 
\bigoplus_{k\in I}\Hom(V_k,V_k),
\quad
(B,i,j)\mapsto
\sum_h\eps(h)B_hB_\oh +ij,$$
where  $\eps$ is any function $\eps\,:\, H\to\CC^\times$ such
that $\eps(h)+ \eps(\oh)=0.$ 
We say that a triple $(B,i,j)\in\mu_{\vb,\wb}^{-1}(0)$ is stable
if there is no nontrivial
$B$-invariant subspace of $\Ker j$.
Let $\mu_{\vb,\wb}^{-1}(0)^s$ be the set of stable triples. 
The group $G_\vb$ acts freely on $\mu_{\vb,\wb}^{-1}(0)^s$.
Put 
$$T(\vb,\wb)=\mu_{\vb,\wb}^{-1}(0)^s/G_\vb,\qquad N(\vb,\wb)=
\mu_{\vb,\wb}^{-1}(0)/\!\!/ G_\vb$$
and let $\pi\,:\, T(\vb,\wb)\to N(\vb,\wb)$ 
be the affinization map 
(it sends $G_\vb\cdot (B,i,j)$ to the only closed
$G_\vb$-orbit contained in $\overline{G_\vb\cdot(B,i,j)}$).
It is proved in [4, 3.10(2)] that 
$T(\vb,\wb)$ is a smooth quasi-projective variety.
Given $\vb^1,\vb^2\in\NN[I]$ 
consider the fiber product
$Z(\vb^1,\vb^2;\wb)=T(\vb^1,\wb)\times_\pi T(\vb^2,\wb).$
Take $\vb^2=\vb^1+\a_k$ where $\a_k$ is a simple root and assume that
$V^1\subset V^2$ have dimension $\vb^1$, $\vb^2$, respectively.
Consider the closed subvariety $C_k^+(\vb^2,\wb)$ of $Z(\vb^1,\vb^2;\wb)$ 
consisting of the pairs of triples $(B^1,i^1,j^1),$ $(B^2,i^2,j^2)$ such that
$B^2_{| V^1}=B^1,\, i^2=i^1,\, j^2_{| V^1}=j^1.$
Put 
$C_k^-(\vb^2,\wb)=
\varphi\bigl(C_k^+(\vb^2,\wb)\bigr)\subset Z(\vb^2,\vb^1;\wb)$ 
where 
$\varphi\,:\,T(\vb^1,\wb)\times T(\vb^2,\wb)\to T(\vb^2,\wb)\times T(\vb^1,\wb)$
permutes the components. The varieties $C_k^\pm(\vb^2,\wb)$ are nonsingular
[4, 5.7].
The group $\tilde G_\wb=G_\wb\times\CC^\times$ acts on $T(\vb,\wb)$ by 
$$(g,t)\cdot (B,i,j)=(tB,t^2ig^{-1},gj),\quad\forall g\in G_\wb,\,
\forall t\in\CC^\times.$$
Let $\Vc_k=\mu_{\vb,\wb}^{-1}(0)^s\times_{G_\vb}V_k$ and 
$\Wc_k$ be respectively the $k$-th tautological bundle 
and the trivial $W_k$-bundle on $T(\vb,\wb).$ 
The bundles $\Vc_k,\Wc_k,$ are $\tilde G_\wb$-equivariant.
The group $\tilde G_\wb$ acts also on $N(\vb,\wb)$,
$C_k^\pm(\vb^2,\wb)$, and $Z(\vb^1,\vb^2;\wb)$.
Let $q$ be the trivial line bundle on $T(\vb,\wb)$ with 
the degree one action of $\CC^\times$.
For any complex $G$-variety X let $K^G(X)$ be the Grothendieck ring of 
$G$-equivariant coherent sheaves on X. Put 
$$\Fc_k(\vb,\wb)=q^{-2}\Wc_k-(1+q^{-2})\Vc_k+q^{-1}
\sum_{h'=k}\Vc_{h''}\in K^{\tilde {G}_\wb}(T(\vb,\wb)).$$
The  rank of $\Fc_k(\vb,\wb)$ is  $(\wb-\vb\,|\,\alpha_k),$ where
$(|)$ is the standard metric on the weight lattice of $\gen.$
We fix a pair of linear maps $\wb\mapsto\wb_\pm$ on the weight 
lattice which are adjoint with respect to $(\quad|\quad)$,
and such that $\wb_++\wb_-=\wb$ for all $\wb$.

\vskip2mm

\noindent{\bf  3. Equivariant  homology and convolution product.} 

\noindent Let $G$ be a complex, connected, linear algebraic group.
For any complex $G$-variety $X$, let $H^G_i(X)$ (resp. $H_G^i(X)$)
be the $i$-th space of 
$G$-equivariant complex Borel-Moore  
homology (resp. of  $G$-equivariant complex cohomology). Put
$$H^G(X)=\bigoplus_i H^G_i(X),\qquad H_G(X)=\bigoplus_i H_G^i(X).$$
See [3] for details on equivariant Borel-Moore homology. 
Let us only recall the following well known facts. 
\vskip1mm

\noindent
- If $Y$ is a closed $G$-subvariety of $X$ and $X$ is smooth,
then $H^G(Y)=H_G(X,X\setminus Y).$ 
Moreover there is a natural map $H_G(X)\to H^G(X).$ 
Call $\a^o\in H^G(X)$ the image of $\a\in H_G(X).$
The $\cup$-product in equivariant cohomology 
induce, via the Poincar\'e duality, a product, noted $\cdot$, in 
equivariant homology. We will denote also by a dot the product
$H_G(X)\otimes H^G(X)\to H^G(X).$
\vskip1mm
\noindent
- Any $G$-equivariant vector bundle $E$ on $X$ admits an equivariant
Chern polynomial $\lambda_z(E)\in H_G(X)[z].$ 
The coefficient of $z$ in $\lambda_z(E)$ is the 
equivariant first Chern class $c_1(E)\in H_G(X).$ The coefficient of 
$z^{\rk(E)}$ in $\lambda_z(E)$ is the equivariant Euler class 
$\lambda(E)\in H_G(X).$
If $E$ is invertible, then $\lambda_z(E)=1+c_1(E)z$. 
Moreover, for any $E$ and $F$
we have $\lambda_z(E\oplus F)=\lambda_z(E)\cup\lambda_z(F)$. The class
$\lambda_z(E)$ depends only on the class of $E$ in $K^G(X).$
\vskip1mm
\noindent
- If $T\subset G$ is a maximal torus, put $\ten=\Lie(T)$. Then  
$H_G^{2i}(pt)=S^{2i}(\ten^*)^W$, where $S^{i}$ is  the 
$i$-symmetric product and $W$ is the Weyl group. 
\vskip1mm
\noindent We will use the following (see [1, Proposition 2.6.47]) :
\vskip1mm
\noindent{\bf Lemma.} {\it Let $X$ be a smooth $G$-variety and let
$C_i$ ($i=1,2$) be two smooth closed $G$-subvarieties. 
Set $C_3=C_1\cap C_2$ and
let $\gamma_i : C_i\hookrightarrow X$ ($i=1,2,3$) be the natural  embedding. 
 Suppose that $C_1$ and
$C_2$ are transversal. 
Then, for all $\a\in H_G(C_1)$ and $\beta\in H_G(C_2)$, 
$$\gamma_{1*}(\a^o)\cdot \gamma_{2*}(\beta^o)=\gamma_{3*}\bigl((\a_{|C_3}\cup \beta_{|C_3})^o
\bigr),$$
where 
$\a_{|C_3}$ (resp. $\beta_{|C_3}$) is
the restriction of $\a$ (resp. $\beta$) to  $C_3$.}\qed 
\vskip1mm
\noindent Let us recall the definition of the convolution product.
Given quasi-projective $G$-varieties $X_1,X_2,X_3,$
consider the projection $p_{ij}\,:\,X_1\times X_2\times X_3\to X_i\times X_j$
for all $1\leq i<j\leq 3$. Consider subvarieties $Z_{ij}\subset X_i\times X_j$
such that the restriction of $p_{13}$ to 
$p_{12}^{-1}Z_{12}\cap p_{23}^{-1}Z_{23}$ is proper and maps to $Z_{13}$. 
The convolution product is the map
$$\star\,:\quad H^G(Z_{12})\otimes H^G(Z_{23})\to H^G(Z_{13}),\quad
\a\otimes \beta\mapsto p_{13\,*}\bigl((p_{12}^*\a)\cdot(p_{23}^*(\beta)\bigr).$$
See [1, 2.7 and the remark (iii), page 113] for more details on
convolution product.
We will essentially consider the case 
$X_i=T(\vb^i,\wb)$ and $Z_{ij}=Z(\vb^i,\vb^j;\wb)$,
where $\vb^1,\vb^2,\vb^3,\wb\in\NN[I]$ and $1\leq i<j\leq 3$.
 
\vskip2mm

\noindent{\bf 4. Statement of the Result.}

\noindent Let $(I,E)$ be a graph of type $ADE$. 
Fix $\wb,\vb^1,\vb^2\in\NN[I]$, with $\vb^2=\vb^1+\a_k$.
For any $k$, denote by $\Vc^1_k$ (resp. $\Vc^2_k$) the vector bundle 
$\Vc_k\boxtimes\Oc_{T(\vb^2,\wb)}$ (resp. $\Oc_{T(\vb^1,\wb)}\boxtimes \Vc_k$)
over $T(\vb^1,\wb)\times T(\vb^2,\wb).$
 The restriction to $C^+_k(\vb^2,\wb)$ of
the sheaf $\Vc^1_k$ is a subsheaf of $\Vc_k^2.$ 
The quotient sheaf $\Lc_k^+=\Vc^2_k/\Vc^1_k$ is 
$\tilde G_\wb$-invariant and invertible. 
Put $\Lc^-_k=\varphi^*\Lc^+_k.$  
Consider the following varieties 
$$N(\wb)=\bigsqcup_\vb N(\vb,\wb),\quad T(\wb)=\bigsqcup_\vb T(\vb,\wb),
\quad Z(\wb)=\bigsqcup_{\vb',\vb''}Z(\vb',\vb'';\wb),$$
where $\vb,\vb',\vb''$ take all the possible values in $\NN[I]$.
Let  $\Delta^\pm$ be  the two natural embeddings 
$$\Delta^+\,:\,C_k^+(\vb^2,\wb)\,\hookrightarrow\, Z(\vb^1,\vb^2;\wb) \and
\Delta^-\,:\,C_k^-(\vb^2,\wb)\,\hookrightarrow\, Z(\vb^2,\vb^1;\wb).$$
If $r\geq 0$, put
$$x_{k,r}^\pm=\sum_{\vb^2}(-1)^{(\a_k\,|\,\vb^2_\pm)}
\Delta^\pm_*(c_1(\Lc_k^\pm)^o)^r\in H^{\tilde G_\wb}(Z(\wb)).\leqno(4.1)$$
Let $\Delta\,:\, T(\vb,\wb)\to T(\vb,\wb)\times T(\vb,\wb)$ be 
the diagonal embedding and set $\hbar=c_1(q^2)^o$. Define
$h_{k,r}$ as the coefficient of 
$\hbar z^{-r-1}$ in 
$$\left( -1+ \sum_\vb\Delta_*{{\lambda_{-1/z}(\Fc_k(\vb,\wb))}
\over{\lambda_{-1/z}(q^2\Fc_k(\vb,\wb))}}\right)^-,\leqno(4.2)$$
where $-$ stands for the expansion at $z=\infty.$
The following result was conjectured by Nakajima 
([5, Introduction], in [2] the result was announced for type $A$). 
\vskip1mm
\noindent{\bf Theorem.} {\it  For all  $\wb\in\NN[I]$, the map 
$\xb^\pm_{k,r}\mapsto x^\pm_{k,r},$ $\hb_{k,r}\mapsto h_{k,r}$
extends uniquely to an algebra homomorphism
$\Phi_\wb\,:\,\Yb_\hbar(\L{\frak g})\to H^{\tilde G_\wb}(Z(\wb)).$}
\qed
\vskip 1mm
\noindent{\bf Remark.} We can prove a similar result for any symmetric
Kac-Moody algebra. Let $A=(a_{kl})_{k,l\in I}$ be a symmetric generalized Cartan matrix.
In the definition of the Yangian, the relation (1.4) becomes
$$\left\{\matrix
[\xb^\pm_{k,r+1},\xb^\pm_{k,s}]-[\xb^\pm_{k,r},\xb^\pm_{k,s+1}]=
\pm\hbar (\xb^\pm_{k,r}\xb^\pm_{k,s}+\xb^\pm_{k,s}\xb^\pm_{k,r})
\hfill\cr\cr
\eta_{-a_{kl}}(z\mp{\hbar\over 2},w)\xb^\pm_k(z)\xb^\pm_l(w)=
\eta_{-a_{kl}}(z,w\mp{\hbar\over 2})\xb^\pm_l(w)\xb^\pm_k(z)
\quad (\text{if}\ k\neq l\hfill)
\endmatrix\right.$$
where
$$
\xb_k^\pm(z)=\sum_{r\geq 0}\xb_{k,r}^\pm\,z^{-r},\and 
\eta_a(z,w)=\prod_{j=1}^n\bigl(z-w+(1+a-2j)\hbar/2\bigr).
$$
In this  case the action of $\CC^\times$ on $T(\vb,\wb)$ and the complex
$\Fc_k(\vb,\wb)$ has to be changed as in [5]. In the proof of the theorem there
are only minor and evident changes to do.

\vskip2mm

\noindent{\bf 5. Proof of the Result.}

\noindent The proof is as in [5, sections 10 and 11] :
we check relations (1.1), (1.2), (1.5) and relations (1.3) and (1.4) in the case
$k\not= l$ by direct computation. Relations (1.3) and (1.4) in the case $k=l$
are proved by reduction to the ${\frak {sl}}_2$-case. We insist here only on the
parts which need different calculations. 

\vskip 2mm

\noindent{\it Relation (1.1).} It is an immediate consequence of the definition,
since for all $x\in H^{\tilde G_\wb}(Z(\vb,\vb';\wb))$ we have
$$h_{k,0}\star x=\rk \Fc_k(\vb,\wb)x=(\wb-\vb|\a_k)x,$$
$$x\star h_{k,0}=\rk \Fc_k(\vb',\wb)x=(\wb-\vb'|\a_k)x.$$

\vskip2mm

\noindent {\it Relation (1.2).} We prove only the plus case,
 the minus being similar.
Fix $\vb^2=\vb^1+
\alpha_l$. We identify $\Fc_k(\vb^1,\wb)$ and $\Fc_k(\vb^2,\wb)$ with their 
pull-back to $C_l^+(\vb^2,\wb)$ via the 1-st and the 2-nd projection. Then,
in $\K^{\tilde G_\wb}(C_l^+(\vb^2,\wb))$, we have 
$$\Fc_k(\vb^1,\wb)-q^2\Fc_k(\vb^1,\wb)
=\Fc_k(\vb^2,\wb)-q^2\Fc(\vb^2,\wb)+[a_{kl}](q^{-1}-q)\Lc_l^+.$$
 It follows that 
$ [h_{k,r},x_{l,s}^+]\in H^{\tilde G_\wb}(C^+_l(\vb^2,\wb))$ 
 is the coefficient of $\hbar z^{-r-1}$ in 
$$\left(\lambda_{-1/z}(\Fc_k(\vb^1,\wb)-q^2\Fc_k(\vb^1,\wb))
-\lambda_{-1/z}(\Fc_k(\vb^2,\wb)-q^2\Fc_k(\vb^2,\wb))\right)^-
x^+_{l,s}=$$
$$=\biggl(\bigl(\lambda_{-1/z}([a_{kl}](q^{-1}-q)\Lc_l^+)-1\bigr)
\lambda_{-1/z}(\Fc_k(\vb^2,\wb)-q^2\Fc_k(\vb^2,\wb))\biggr)^-
x^+_{l,s}.$$
Set $$A_s=
\lambda_{-1/z}(\Fc_k(\vb^2,\wb)-q^2\Fc_k(\vb^2,\wb))
x^+_{l,s},$$
$$
X=\lambda_{-1/z}([a_{kl}](q^{-1}-q)\Lc_l^+)={{1-(c_l^+-a_{kl}\hbar/2)z^{-1}}\over{
1-(c_l^++a_{kl}\hbar/2)z^{-1}}}.$$
Then the LHS and the RHS of the relation (1.2) are respectively equal to the coefficient of 
$\hbar z^{-r-1}$ in
$$\left(2z(X-1)A_s-2(X-1)A_{s+1}\right)^-=\left(2(X-1)(z-c_l^+)A_s\right)^-
 \and \left(\hbar a_{kl}(X+1)A_s\right)^-.$$
We are then reduced to the identity, easily checked, in $H^{\tilde G_\wb}(C_l^+(\vb^2,\wb))$ :
$$2(X-1)(z-c_l^+)=\hbar a_{kl}(X+1).$$

\vskip2mm

\noindent{\it Relation (1.3) with $k\not= l.$}  
Fix $\vb^1,\vb^2,\vbt^2,\vb^3$, such
that
$$\vbt^2=\vb^1-\a_l=\vb^3-\a_k=\vb^2-\a_k-\a_l.$$
If $1\leq i<j\leq 3$, consider the projections
$$ p_{ij}\,:\, T(\vb^1,\wb)\times T(\vb^2,\wb)\times T(\vb^3,\wb)\to
T(\vb^i,\wb)\times T(\vb^j,\wb),$$
$$\tilde p_{ij}\,:\, T(\vb^1,\wb)\times T(\vbt^2,\wb)\times T(\vb^3,\wb)\to
T(\vbt^i,\wb)\times T(\vbt^j,\wb),$$
where we set $\vbt^1=\vb^1,\vbt^3=\vb^3.$
 We have
$$x_{k,r}^+\star x_{l,s}^-=(-1)^{(\a_k|\vb_+^2)+(\a_l|\vb^3_-)}
p_{13*}\left(p_{12}^*(c_1(\Vc^2_k/\Vc^1_k)^{o r})\cdot 
p_{23}^*(c_1(\Vc^3_l/\Vc_l^2)^{o s})\right),$$
$$x_{l,s}^-\star x_{k,r}^+=(-1)^{(\a_l|\vbt^2_-)+(\a_k|\vb_+^3)}
\tilde p_{13*}\left(\tilde p_{12}^*(c_1(\tilde\Vc^2_l/\Vc_l^1)^{o s})\cdot \tilde
p_{23}^*(c_1(\Vc_k^3/\tilde\Vc_k^2)^{o r})\right).$$
It is proved in [5, Lemma 10.2.1] that  the intersections 
$$p_{12}^{-1}C_k^+(\vb^2,\wb)\cap p_{23}^{-1}C_l^-(\vb^3,\wb)\and
\tilde p_{12}^{-1}C_l^-(\vbt^2,\wb)\cap\tilde p_{23}^{-1}C_k^+(\vb^3,\wb)$$
are transversal in $T(\vb^1,\wb)\times T(\vb^2,\wb)\times T(\vb^3,\wb)$ and
$T(\vb^1,\wb)\times T(\vbt^2,\wb)\times T(\vb^3,\wb)$ 
and that there exists a $\tilde G_\wb$-equivariant isomorphisms between them 
which induces the isomorphisms : 
$$ \Vc^2_k/\Vc^1_k\simeq \Vc^3_k/\tilde\Vc^2_k\and
\tilde\Vc^2_l/\Vc^1_l\simeq \Vc^3_l/\Vc^2_l.$$
The result follows from the lemma in section 3.

\vskip2mm

\noindent{\it Relation (1.4) with $k\not=l$.} We prove only the plus case.
Fix $\vb^1,\vb^2,\vbt^2,\vb^3,$
such that
$$\vb^3=\vbt^2+\a_k=\vb^2+\a_l=\vb^1+\a_k+\a_l.$$
Consider the projections $p_{ij}$ and $\tilde p_{ij}$ ($1\leq i<j\leq 3$) as before. 
The intersections 
$$Z_{kl}=p_{12}^{-1}C_k^+(\vb^2,\wb)\cap p_{23}^{-1}C_l^-(\vb^3,\wb)\and
Z_{lk}=\tilde p_{12}^{-1}C_l^-(\vbt^2,\wb)\cap\tilde p_{23}^{-1}C_k^+(\vb^3,\wb)$$
are transversal in $T(\vb^1,\wb)\times T(\vb^2,\wb)\times T(\vb^3,\wb)$ and
$T(\vb^1,\wb)\times T(\vbt^2,\wb)\times T(\vb^3,\wb)$ 
(see [5, Lemma 10.3.1]). Since $k\not= l$, the restriction of $p_{13}$
and $\tilde p_{13}$ to $Z_{kl}$ and $Z_{lk}$ is an embedding into 
$Z(\vb^1,\vb^3;\wb).$ Call it $\iota_{kl}$ and $\iota_{lk}$ respectively.
Put  $b_k=c_1(\Vc_k^3-\Vc_k^1)$, $b_l=c_1(\Vc_l^3-\Vc_l^1).$
 We have (see the lemma in section 3)
$$x_{k,r}^+\star x_{l,s}^+=(-1)^{(\a_k|\vb_+^2)+(\a_l|\vb^3_+)}
\iota_{kl*}\bigl((p_{12}^*(b_k^r)_{|Z_{kl}}\cup 
p_{23}^*(b_l^s)
_{|Z_{kl}})^o\bigr),$$
$$x_{l,s}^+\star x_{k,r}^+=(-1)^{(\a_l|\vbt^2_+)+(\a_k|\vb_+^3)}
\iota_{lk*}\bigl((p_{12}^*(b_k^r)_{|Z_{lk}}\cup 
p_{23}^*(b_l^s)
_{|Z_{kl}})^o\bigr).$$
Take $h\in H$ such that $h'=l$ and $h''=k.$ The map $B_\oh$ 
may be viewed as a section of the $\tilde G_\wb$-bundle $\Ec_{kl}=q(\Vc^3_l/\Vc_l^1)^*
\otimes (\Vc_k^3/\Vc_k^1)$ on $p_{13}(Z_{kl})$ (where we set $\Ec_{kl}=0$ if $a_{kl}=0$). 
Similarly 
$B_h$ is a section of
the $\tilde G_\wb$-bundle $\Ec_{lk}=q(\Vc^3_k/\Vc_k^1)^*
\otimes (\Vc_l^3/\Vc_l^1)$ on $\tilde p_{13}(Z_{lk})$ (where again we set $\Ec_{lk}=0$ if
$a_{kl}=0$).
In [5, 10.3.9] it is proved that $B_\oh$ and $B_h$ are transversal to the 
zero section respectively. Moreover 
$$p_{13}(Z_{kl})\cap B_\oh^{-1}(0)=
\tilde p_{13}(Z_{lk})\cap B_h^{-1}(0).$$
Then,  
$$\iota_{kl*}(c_1(\Ec_{kl})^o)\, x_{k,r}^+\star x_{l,s}^+=(-1)^{a_{kl}}
\iota_{lk*}(c_1(\Ec_{lk})^o)\,x_{l,s}^+\star x_{k,r}^+ ,$$
i.e.
$$\iota_{kl*}(b_k^o-b_l^o+\hbar/2)\,x_{k,r}^+\star x_{l,s}^+ =
\iota_{lk*}(b_l^o-b_k^o+\hbar/2)\,x_{l,s}^+\star x_{k,r}^+ .$$
The relation (1.4) follows immediately from this.

\vskip2mm

\noindent{\it Relations (1.3) and (1.4) with $k=l$.} Thank to the same argument
than in [5, 11.3] we are reduce to the case of 
$(I,E)$ of type $A_1.$ In this case $\vb$ and $\wb$
are identified with natural numbers, so let us call them $v$ and $w$.
Moreover we will omit everywhere the subindex 1.
Let $\Gr_v(w)$ be the variety of $v$-dimensional subspaces in $W$.
It is easy to see that $T(v,w)\simeq T^*\Gr_v(w)$. 
The group $G_w$ acts in the obvious way on $T(v,w)$.
The group $\CC^\times$ acts by scalar multiplication on the fibers
of the cotangent bundle. Fix $T_1,...,T_w$ such that
$$\matrix
K^{\tilde G_w}(\Gr_v(w))=\CC[q^{\pm 1},T_1^{\pm 1},...,T_w^{\pm 1}]^{S_v\times S_{w-v}},
\hfill\cr\cr
\wedge^i\Vc=e_i(T_1,...,T_v),\qquad\wedge^i\Wc=e_i(T_1,...,T_w),
\hfill\endmatrix$$
where $e_i$ is the $i$-th elementary symmetric function. 
We get
$$\Fc(v,w)=q^{-2}\Wc-(1+q^{-2})\Vc=q^{-2}(T_{v+1}+\cdots
+T_w)-(T_1+\cdots+ T_v).$$
Put $t_k=c_1(T_k)^o$. Then $H^{\tilde G_w}(T(w))=
\oplus_{v=0}^w\CC[\hbar,t_1,\ldots t_w]^{S_v\times S_{w-v}}.$ The following lemma 
is proved as in [1, Claim 7.6.7].
\vskip1mm
\noindent{\bf Lemma.}{\it The space $H^{\tilde G_w}(T(w))$ is a faithful
module over $H^{\tilde G_w}(Z(w))$. }\qed
\vskip1mm
\noindent 
The operators 
$x_r^\pm$ on $H^{\tilde G_w}(T(w))$ can be written down
explicitely. Put
$$O(v,w)=\{(V^1,V^2)\in\Gr_{v-1}(w)\times\Gr_v(w)\,|\, 
V^1\subset V^2\}.$$
The Hecke correspondence $C^+(v,w)$ is the conormal bundle 
to $O(v,w)$. Consider the projections $p_1,p_2$ from $O(v,w)$
to the first and the second component and let 
$\pi\,:\, T(v,w)\to Gr_v(w)$ be the projection.
We can prove as in [6, Lemme 5] that if $\a\in H_{\tilde G_w}(O(v,w))$ and 
$\beta\in H^{\tilde G_w}(\Gr_v(w))$, then
$$\pi^*(\a)\star\pi^*(\beta)=
p_{1 *}(\lambda(q^2T^*p_1)\cdot \a\cdot{p_2}^*\beta),$$ 
where $T^*p_1$ is the relative cotangent bundle to $p_1$.
The map $p_1$ is a $\PP^{w-v}$-fibration, then we have 
$$\lambda(q^2T^*p_1)=\prod_{m=v+1}^w(t_m-t_v+\hbar)
\in H^{\tilde G_w}(O(v,w)).$$
Let us introduce the following notation. Fix $z\in [1,w]=\{1,2,...,w\}$ and let 
$I=(I_1,I_2)$ be a partition of $[1,w]$ into two subset of cardinality $z$ and $w-z$
respectively, say $I_1=\{a_1,a_2,...,a_z\},$ $I_2=\{b_1,b_2,...,b_{w-z}\}.$ Then put
$$f(t_{I_1};t_{I_2})=f(t_{a_1},t_{a_2},...,t_{a_z},t_{b_1},t_{b_2},...,t_{b_{w-z}}).$$
Thus (see [6, Lemme 1]), for any 
$f\in\CC[\hbar,t_1,...,t_w]^{S_v\times S_{w-v}},$
\vskip2mm
\noindent
$$x^+_r(f)(t_{[1,v-1]};t_{[v,w]})= \sum_{k=v}^w
f(t_{[1,v-1]\cup\{k\}};t_{[v,w]\setminus\{k\}})t_k^r
\prod_{m\in[v,w]\setminus\{k\}}\left(1+{\hbar\over {t_k-t_m}}\right),\leqno(5.1)$$
$$x^-_r(f)(t_{[1,v+1]};t_{[v+2,w]})= \sum_{k=1}^{v+1}
f(t_{[1,v+1]\setminus\{k\}};t_{[v+2,w]\cup\{k\}})t_k^r
\prod_{m\in [1,v+1]\setminus\{k\}}\left(1+{\hbar\over {t_m-t_k}}\right).\leqno(5.2)$$
We have
$$\lambda_{-z}\Fc(v,w)=
{{\prod_{m=v+1}^w\bigl(1-z(t_m-\hbar)\bigr)}\over
{\prod_{m=1}^v\bigl(1-z(t_m+1)\bigr)}}.$$
Thus $h_r(f)$ is the coefficient of $\hbar z^{-r-1}$ in 
$$f\left(\prod_{m=1}^v{{z-t_m-\hbar}\over{z-t_m}}\prod_{m=v+1}^w
{{z-t_m+\hbar}\over{z-t_m}}\right)^-.\leqno(5.3)$$
\vskip1mm
\noindent{\bf Proposition.} {\it Relations $(1.3)$ and $(1.4)$ hold in the ${\frak{sl}}_2$-case. }

\vskip1mm

\noindent{\it Proof.} 
Let us prove relation (1.3).
Fix $f\in\CC[\hbar,t_1,...,t_w]^{S_v\times S_{w-v}}.$
Using formulas (5.1) and (5.2), we have 
$$\big(x^-_sx^+_r(f)\big)
(t_{[1,v]};t_{[v+1,w]})= \sum_{l=1}^v\sum_{k\in [v+1,w]\cup\{l\}}
f(t_{([1,v]\setminus\{l\})\cup\{k\}};t_{([v+1,w]\cup\{l\})\setminus\{k\}})
t_l^st_k^rX_{kl},$$
$$\big(x^+_rx^-_s(f)\big)
(t_{[1,v]};t_{[v+1,w]})= \sum_{l\in[1,v]\cup\{k\}}\sum_{k=v+1}^w
f(t_{([1,v]\cup\{k\})\setminus\{l\}};t_{([v+1,w]\setminus\{k\})\cup\{l\}})
t_l^st_k^rY_{kl},$$
where
$$\matrix
{\ds X_{kl}=
\prod_{m\in [1,v]\setminus\{l\}}\left(1+{\hbar\over {t_m-t_l}}\right)
\prod_{n\in[v+1,w]\cup\{l\}\setminus\{k\}}\left(1+{\hbar\over
{t_k-t_n}}\right)},\hfill\cr\cr
{\ds Y_{kl}=\prod_{m\in[1,v]\cup\{k\}\setminus\{l\}}
\left(1+{\hbar\over {t_m-t_l}}\right)
\prod_{n\in[v+1,w]\setminus\{k\}}\left(1+{\hbar\over{t_k-t_n}}\right).}
\endmatrix$$
The terms with $k\not= l$ cancel out in the bracket. We get
$$[x^+_r,x^-_s](f)= f \sum_{k=v+1}^wt_k^st_k^r
\prod_{m\in[1,v]}\left(1+{\hbar\over{t_m-t_k}}\right)
\prod_{n\in[v+1,w]\setminus\{k\}}\left(1+{\hbar\over{t_k-t_n}}\right)-$$
$$-f \sum_{l=1}^vt_l^st_l^r
\prod_{m\in[1,v]\setminus\{l\}}\left(1+{\hbar\over{t_m-t_l}}\right)
\prod_{n\in[v+1,w]}\left(1+{\hbar\over{t_l-t_n}}\right).$$
Put 
$$A(z)=\prod_{m=1}^w(z-t_m),\qquad
B(z)=\prod_{m=1}^v(z-t_m-\hbar)\prod_{m=v+1}^w(z-t_m+\hbar).$$
Then it is easy to check that 
$$\hbar[x^+_r,x^-_s](f)=f\sum_{k=1}^w
t_k^{r+s}{{B(t_k)}\over {A'(t_k)}}=
f\Res_\infty z^{r+s}{{B(z)}\over{A(z)}}.$$
This is the definition of $h_{r+s}(f)$ given in (5.3).
As for the relation (1.4), note that,
using (5.1), we get 
\vskip1mm
\noindent$\big(x^+_sx^+_r(f)\big)
(t_{[1,v-2]};t_{[v-1,w]})=$\hfill\break
$$=\sum_{l=v-1}^w\sum_{k\in[v-1,w]\setminus\{l\}}
f(t_{[1,v-2]\cup\{k,l\}};t_{[v-1,w]\setminus\{k,l\}})
t_l^st_k^rZ_{kl},$$
where
$$Z_{kl}=\prod_{n\in[v-1,w]\setminus\{l\}}\left(1+{\hbar\over{t_l-t_n}}\right)
\prod_{m\in[v-1,w]\setminus\{k,l\}}\left(1+{\hbar\over{t_k-t_m}}\right).$$
The relation  in the plus case follows now by a direct computation. 
\qed

\vskip2mm

\noindent{\it Relation (1.5).} The proof is exactly as in [5, 10.4], so we omit it .

\vskip1mm
\noindent{\bf Remark.} Nakajima [5, Theorem 9.4.1] has proved that there
exists an algebra morphism
$$\Psi_\wb\,:\,\Ub_q(\L\gen)\to K^{\tilde G_\wb}(Z(\wb))\otimes_{\ZZ[q,q^{-1}]}\CC(q),$$
where the algebra to the left is the quantized enveloping algebra of $\L\gen$ and
the algebra to the right is equipped with the convolution product. Using $\Phi_\wb$
and $\Psi_\wb$ we can construct the finite dimensional simple modules of $\Yb_\hbar(\L\gen)$
and $\Ub_q(\L\gen)$ respectively (see [5, section 14]).
In particular $\Yb_\hbar(\L\gen)$
and $\Ub_q(\L\gen)$ have the same finite dimensional representation theory.
More precisely let ${\frak C}$ (resp. ${\frak D}$) be the abelian category
of finite dimensional $\Ub_q(\L\gen)$-(resp. $\Yb_\hbar(\L\gen)$-)modules
such that the Drinfeld polynomials of the simple factors
 have roots in $q^\ZZ$ (resp.  $\ZZ$). 
\vskip1mm
\noindent{\bf Proposition.} {\it  The characters (as $\Ub_q(\gen)$-modules
and $\Ub(\gen)$-modules resp.)
of the simple finite dimensional modules in ${\frak C}$ and in
${\frak D}$
are the same.} \qed

\vskip2mm

\noindent{\bf Aknowledgments.}
The author thanks E. Vasserot for useful discussions.

\vskip2mm

\noindent{\bf References.}

\noindent 1. Chriss, N., Ginzburg, V.: {\sl Representation theory and complex
geometry}, Birkh\"auser, Boston-Basel-Berlin, 1997.

\vskip 1mm

\noindent 2. Ginzburg, V., Vasserot, E.: Langlands reciprocity for 
affine quantum groups of type $A_n$, {\sl Internat. Math. Res. Notices},
{\bf 3}, 1993, 67-85.

\vskip1mm

\noindent 3. Lusztig, G.: Cuspidal local systems and graded Hecke 
algebras I, {\sl Inst. Hautes \'Etudes Sci. Publ. Math.}, {\bf 67},
1988, 145-202.

\vskip1mm

\noindent 4. Nakajima, H.: Quiver varieties and Kac-Moody algebras,
{\sl Duke. Math. J.}, {\bf 91}, 1998, 515-560.

\vskip1mm

\noindent 5. Nakajima, H.: Quiver varieties and finite dimensional 
representations of quantum affine algebras,  Preprint QA/9912158.

\vskip1mm

\noindent 6. Vasserot, E.: Repr\'esentations de groupes quantiques
et permutations, {\sl  Ann. Scient. \'Ec. Norm. Sup.}, {\bf 26}, 1993, 747-773.

\enddocument